\documentclass[12pt]{smfart}

\usepackage[frenchb]{babel}
\usepackage{mathptmx}
\usepackage{amsmath}
\usepackage{amssymb}
\usepackage{smfthm}
\usepackage{amscd}
\usepackage{smfenum}

\usepackage{amssymb,amsfonts}

\usepackage{diagrams}
\usepackage{graphics}
\usepackage[T1]{fontenc}     
\usepackage[latin1]{inputenc}

\newcounter{ctrmasection}
\newcounter{ctrmasoussection}

\newtheorem{theoreme}{Th{\'e}or{\`e}me}[section]
\newtheorem{lemme}[theoreme]{Lemme}

\newenvironment{pf*}{  \bf{D{\'e}monstration : }\rm}{$\Box$}

{\begin{displaymath} \begin{array}{c}}%
{\end{array} \end{displaymath}}
\newcommand{\Aut}{\operatorname{Aut}}

\newcommand{\Pic}{\operatorname{Pic}}

\newcommand{\B}{\operatorname{B}}
\newcommand{\Cplus}{\operatorname{C_+}}

\begin{document}
\title[Points rationnels et automorphismes affines]{G{\'e}om{\'e}trie, points rationnels et it{\'e}r{\'e}s des automorphismes de l'espace affine} 
\alttitle{Geometry, rational points, and iterates of automorphisms of the affine space}
\author{Sandra Marcello}

\address{Max-Planck-Institut f{\"u}r Mathematik,
Vivatsgasse 7, 53111, Bonn, Deutschland \\}
\email{marcello@mpim-bonn.mpg.de}

\date{\today}

\begin{abstract}Nous {\'e}tudions  sur un corps de nombres les it{\'e}r{\'e}s
  des automorphismes de l'espace affine. Nous nous int{\'e}ressons plus pr{\'e}cis{\'e}ment aux points p{\'e}riodiques et aux points non-p{\'e}riodiques; pour les premiers les questions sont analogues {\`a} celles pos{\'e}es pour les points de torsion des vari{\'e}t{\'e}s ab{\'e}liennes, pour les seconds les questions sont comparables aux probl{\`e}mes de d{\'e}compte de points rationnels sur les vari{\'e}t{\'e}s. 
Dans le cadre de cette
  {\'e}tude apparaissent de mani{\`e}re naturelle des invariants
  g{\'e}om{\'e}triques qui sont d{\'e}finis {\`a} partir  des cones de diviseurs de vari{\'e}t{\'e}s construites {\`a} partir des automorphismes de l'espace affine. Nous {\'e}tudions et majorons en toute dimension ces invariants, nous obtenons pour l'un d'entre eux une borne qui d{\'e}pend du degr{\'e} dynamique, borne qui est atteinte pour certains automorphismes, en particulier la majoration est optimale. Pour un autre des invariants, nous obtenons en dimension quelconque gr{\^a}ce {\`a} une construction g{\'e}om{\'e}trique et sous certaines conditions, une meilleure borne. 
\end{abstract}

\begin{altabstract} We study  over a number field, the iterates of
  automorphisms of the affine space. More precisely, we are interested in the periodic and non-periodic points; for the former the questions are similar to the ones about torsion points on abelian varieties, for the latter the questions are similar to the problems on counting rational points on varieties. 
 In order to study this problem,
  we define geometric invariants, which are defined thanks to the cone of divisors  associated to some varieties constructed from the automorphisms. We study, in all dimension, these invariants and bound them above for one of them we obtain a bound depending  on the dynamical degree, this bound is reach by some automorphisms, in particular this bound is optimal. For an other invariant, we obtain in all dimension, thanks to a geometrical construction and under some geometric conditions, a better bound.

\end{altabstract}


\altkeywords{Heights,\  iterates,automorphisms of affine space,\ dynamical degree,\  ample cone, \ effectice cone}

\maketitle

\tableofcontents

\mainmatter

\section{Introduction}
Les conjectures maintenant classiques de Manin (voir \cite{BM} et \cite{fmt}) proposent un lien entre la g{\'e}om{\'e}trie d'une vari{\'e}t{\'e} alg{\'e}brique et une estimation asymptotique du nombre de ses points rationnels. Pour une synth{\`e}se et les derniers d{\'e}veloppements voir le s{\'e}minaire Bourbaki \cite{Peyre}.
\noindent Ainsi, dans \cite{BM}, dans le cadre de probl{\`e}mes autour du d{\'e}compte de points rationnels sur des vari{\'e}t{\'e}s, V. Batyrev and Y. Manin
d{\'e}finissent un invariant arithm{\'e}tique et un invariant g{\'e}om{\'e}trique ce dernier co{\"\i}ncide dans le cas des diviseurs amples avec l'invariant de Nevanlinna.
Dans ce m{\^e}me texte \cite{BM}, les auteurs calculent ces invariants dans certains cas  et {\'e}noncent des conjectures les concernant; ainsi que nous  ayons {\'e}t{\'e} amen{\'e}e {\`a} d{\'e}finir des invariants g{\'e}om{\'e}triques des automorphismes de l'espace affine dans le cadre de l'{\'e}tude des points rationnels p{\'e}riodiques  ou non de ces applications n'a rien de surprenant, comme nous le montre le th{\'e}or{\`e}me A.
Pour plus de pr{\'e}cisions sur l'invariant de Nevanlinna nous renvoyons le lecteur {\`a} \cite{hs}  F.5.4.

\noindent Dans ce texte, nous nous int{\'e}ressons notamment {\`a} une question du m{\^e}me type dans le cadre des it{\'e}r{\'e}s d'automorphismes de l'espace affine. Plus pr{\'e}cis{\'e}ment nous nous int{\'e}ressons {\`a} un probl{\`e}me de d{\'e}compte de points rationnels dans les orbites des points non-p{\'e}riodiques et  cherchons {\`a} trouver un lien avec la g{\'e}om{\'e}trie des applications que nous consid{\'e}rons. Nous nous int{\'e}ressons {\'e}galement {\`a} la conjecture de Silverman sur les points p{\'e}riodiques isol{\'e}s.

\begin{defi}Soit $k$ un corps de nombres. Soit $\phi\in\Aut(\mathbb{A}^r(k))$.
 Soit $P\in \mathbb{A}^r(k)$. Ce point $P$  
est un {\it{point $\phi$-p{\'e}riodique isol{\'e}}}
si et seulement s'il existe $n \in \mathbb{N}\setminus \{0\}$ tel que
$P$ est isol{\'e} (au sens de la topologie de Zariski) dans
 $$\{Q\in \mathbb{A}^r(k)\quad \mbox{tel que}\quad \phi ^{n}(Q)=Q\}.$$ 
\end{defi}  

\noindent{\textbf{Conjecture de Silverman:}} {\it{Soit $k$ un corps de nombres. Soit $\phi\in Aut(\mathbb{A}^r(k))$.
L'ensemble des points p{\'e}riodiques isol{\'e}s de $\phi$ est un ensemble fini.}}

\vspace{0.5cm}

Pour une formulation plus motiv{\'e}e des questions que nous consid{\'e}rons, nous revoyons le lecteur {\`a} \cite{1art}.
Dans \cite{2art}, nous consid{\'e}rons une question analogue dans le cas des corps $p$-adiques, question que nous relions {\`a} la conjecture de Silverman.
Nous introduisons maintenant les notations et objets n{\'e}cessaires notamment {\`a} la construction des invariants g{\'e}om{\'e}triques, puis donnons les r{\'e}sultats reliant arithm{\'e}tique et dynamique.
   
\noindent Nous nous pla{{\c c}}ons pour les aspects g{\'e}om{\'e}triques et arithm{\'e}tiques  respectivement sur $\mathbb{C}$ et $\overline{\mathbb{Q}}$.

Nous notons $\mathbb{A}^r$ (resp. $\mathbb{P}^r$) l'espace affine (resp. projectif) de
dimension $r\geq 2$ et
$H$ l'hyperplan {\`a} l'infini. Lorsque $\phi\in\Aut
(\mathbb{A}^r)$ est un automorphisme de $\mathbb{A}^r$, nous noterons encore
$\phi:\mathbb{P}^r\dots\rightarrow \mathbb{P}^r$ l'application rationnelle induite et d{\'e}signerons par
$Z(\phi)$ le {\it{lieu de non-d{\'e}finition}}, {\it{i.e}} le lieu g{\'e}om{\'e}trique, contenu dans $H$, o{\`u} $\phi$ n'est pas d{\'e}finie.

\noindent De plus, si $\phi\in\Aut(\mathbb{A}^r)$, alors $\phi$ est
d{\'e}fini par des polyn{\^o}mes $P_1,\dots,P_r$. On d{\'e}finit le {\it{degr{\'e}
    alg{\'e}brique}} ou  {\it{degr{\'e}}} de $\phi$ comme suit : $deg(\phi)=d=\max_i(d_i)$ avec $d_i$ degr{\'e} total de $P_i$.

\begin{defi}Soit V une vari{\'e}t{\'e} lisse projective. On dit que le morphisme birationnel $\pi:V\rightarrow\mathbb{P}^r$  est une {\it
r{\'e}solution} de $\phi$ (not{\'e}e  $(\pi,V)$) si $\phi\circ\pi$ et
$\phi^{-1}\circ\pi$ sont des morphismes de $V$ vers $\mathbb{P}^r$.
\end{defi}

\begin{rema}
Parfois nous noterons:
$\psi=\phi\pi$ et $\psi'=\phi^{-1}\pi$.
\end{rema} 

Une telle r{\'e}solution existe toujours comme nous le voyons, en prenant pour $V$
la d{\'e}singularisation de 
l'adh{\'e}rence de Zariski dans $(\mathbb{P}^r)^3$ de
$$\{(x,\phi(x),\phi^{-1}(x))\;|\;x\in\mathbb{A}^r\}$$
et pour $\pi=pr_1$ puisqu'alors $\phi\circ\pi=pr_2$ et
$\phi^{-1}\circ\pi=pr_3$.

\paragraph*{Notations} Si $V$ est une vari{\'e}t{\'e}, on note respectivement $\Pic_{\mathbb{Q}}(V)$ et $\Pic_{\mathbb{R}}(V)$ son
groupe de Picard tensoris{\'e} respectivement par $\mathbb{Q}$  et $\mathbb{R}$, $\Pic^+_{\mathbb{Q}}(V)$ et $\Pic^+_{\mathbb{R}}(V)$ les cone ferm{\'e}s
engendr{\'e}s par les classes de diviseurs effectifs, $\Pic^{a}_{\mathbb{R}}(V)$ le cone (ouvert)
engendr{\'e} par les classes de diviseurs amples et enfin $\Pic^{nef}_{\mathbb{R}}$ le cone engendr{\'e} par les classes de diviseurs nef (num{\'e}riquement effectif).

Soit $(\pi,V)$ une r{\'e}solution de $\phi$ et  $\alpha\in\mathbb{R}$, nous notons:
$$D(\alpha,\pi,V):=(\phi\circ\pi)^*H+(\phi^{-1}\circ\pi)^*H-\alpha  \pi^*H.$$

Pour des raisons arithm{\'e}tiques qui apparaissent clairement dans l'{\'e}nonc{\'e} du th{\'e}or{\`e}me B, nous introduisons :

$$\alpha_{max,eff}(\phi,\pi,V):=\sup\{\alpha\in\mathbb{R}\;|\;D(\alpha,\pi,V) \in
\Pic^+_{\mathbb{R}}(V)\}, $$
$$\alpha_{max,amp}(\phi,\pi,V):=\sup\{\alpha\in\mathbb{R}\;|\;D(\alpha,\pi,V) \in
\Pic^{a}_{\mathbb{R}}(V)\}.$$
$$\alpha_{max,nef}(\phi,\pi,V):=\sup\{\alpha\in\mathbb{R}\;|\;D(\alpha,\pi,V) \in
\Pic^{nef}_{\mathbb{R}}(V)\}.$$

Nous d{\'e}finissons {\'e}galement :
$$\alpha(\phi,eff):=\sup\{\alpha_{\max, eff}(\phi,\pi,V)| \pi, V\},$$
$$\alpha(\phi,amp):=\sup\{\alpha_{\max,amp}(\phi,\pi,V)| \pi, V\},$$
$$\alpha(\phi,nef):=\sup\{\alpha_{\max,nef}(\phi,\pi,V)| \pi, V\}.$$

Nous montrons que ces nombres ne d{\'e}pendent pas de la r{\'e}solution choisie, ce sont donc des invariants g{\'e}om{\'e}triques des automorphismes de l'espace affine.
Nous les appelerons respectivement l'{\it{indice effectif}}, l'{\it{indice ample}} et l'{\it{indice nef}} d'un automorphisme de l'espace affine.
Nous montrons que les indices ample et nef sont {\'e}gaux.

\noindent Nous utiliserons l'indice ample pour des applications arithm{\'e}tiques, dans l'optique de ces m{\^e}mes applications arithm{\'e}tiques, nous d{\'e}finissons: {\it{l'indice effectif arithm{\'e}tique}}.

\noindent Soit $D\in\Pic(V)$, $m\in\mathbb{N}\setminus\{0\}$.
 Nous noterons $\mid D\mid$ le syst{\`e}me lin{\'e}aire associ{\'e} au diviseur $D$  et $\B_m(D)=\cap_{\Delta\in\mid mD\mid} \Delta$ les points base du syst{\`e}me lin{\'e}aire $\mid mD\mid $. Nous d{\'e}finissons : 
$$\B_{\infty}(D)=\cap_{m\geq 1,mD\in \Pic(V)}\B(mD).$$ 
Nous notons,
$$\Cplus=\{D\in \Pic^+_{\mathbb{Q}}\quad \mid \quad \B_{\infty}(D)\subset \pi^{-1}(H)\}.$$

\begin{rema}
La propri{\'e}t{\'e} $\B_{\infty}(D)\subset \pi^{-1}(H)$ est {\'e}quivalente {\`a} :
\begin{center} 
il existe $m\geq 1$ tel que $\B_{m}(D)\subset \pi^{-1}(H)$.
\end{center}
\end{rema}
Pour nos applications arithm{\'e}tiques, nous nous int{\'e}ressons {\`a} :
$$D(\alpha,\pi,V):=(\phi\circ\pi)^*H+(\phi^{-1}\circ\pi)^*H-\alpha  \pi^*H\in \Cplus.$$
Si pour tout $\alpha>0$ cette propri{\'e}t{\'e} n'est pas v{\'e}rifi{\'e}e nous posons alors  $\alpha_{arith}(\phi,\pi,V)=0$, dans le cas contraire nous d{\'e}finissons 
$$\alpha_{arith}(\phi,\pi,V)=\sup\{\alpha\in\mathbb{Q}\quad \mid\quad D(\alpha,\pi,V)\in \Cplus \}.$$

Nous appellerons ce nombre, qui d{\'e}pend {\it{a priori}} de la r{\'e}solution choisie, l'{\it{indice effectif arithm{\'e}tique}} de $\phi$ associ{\'e} {\`a} la r{\'e}solution $(\pi,V)$.

La lettre $h$ d{\'e}signe la hauteur (de Weil) logarithmique normalis{\'e}e usuelle
(voir par exemple \cite{hs}).

\paragraph*{Notations} Nous d{\'e}finissons une fonction de
comptage $N(\phi,P,B)$ pour les orbites des points non-p{\'e}riodiques.
$$N(\phi,P,B)=N(B)=\#\{\phi^n(P)\mid n \in \mathbb{Z},\ h(\phi^n(P))\leq B\}.$$
Pour $B$ grand, nous dirons que $N(B)$ est {\it{comparable {\`a}}} $\log B$ s'il existe 

\noindent $a_1,a_2,a_3,a_4\in \mathbb{R}$ tels que\,: $$a_1\log B+a_2\leq N(B)\leq a_3\log B +a_4.$$
Si $\lim \frac{N(B)}{\log B}=a$ alors, $N(B)$ est dit {\it{{\'e}quivalent {\`a}}} $\log B$.
On les note respectivement $N(B)\gg\ll \log B$ et $N(B)\sim a \log B$.

\begin{defi}
Soit $P=(x_0,\cdots,x_r)\in\mathbb{P}^r(\overline{\mathbb{Q}})$,on appelle {\it{degr{\'e} }} de $P$, le degr{\'e} de $\mathbb{Q}(P)$   \it{i.e}$[\mathbb{Q}(P):\mathbb{Q}]$, o{\`u} 
$ \mathbb{Q}(P)=\mathbb{Q}\big(\frac{x_i}{x_j},\ \mbox{avec}\ 0\leq i,j\leq r \ \mbox{et}\ x_j\neq 0\big)$.
\end{defi}

\paragraph*{{\bf{Th{\'e}or{\`e}me A}}} \label{eclat1}
  {\it{Soit $\phi$ un automorphisme de ${\mathbb{A}}^r$ de degr{\'e}
  $d\geq 2$.
Supposons v{\'e}rifi{\'e}e l'une des trois conditions suivantes.
\begin{itemize}
\item[(H1)] Il existe une r{\'e}solution  $(\pi,V)$ de $\phi$ pour laquelle : $\alpha_{arith}(\phi,\pi,V)>2$.
\item[(H2)] Soit $\alpha(\phi,amp)>2$.
\item[(H3)] Soit $\alpha(\phi,nef)>2$.
 \end{itemize}
\noindent Alors~}}:
\begin{enumerate}
\item[(P1)] {\it{l'ensemble des points p{\'e}riodiques de $\phi$ est de hauteur born{\'e}e sur $\overline{\mathbb{Q}}$, {\it{a fortiori}} l'ensemble des points p{\'e}riodiques de $\phi$ de hauteur et de degr{\'e} born{\'e}s est fini;}}
\item[(P2)] {\it{si P est non-p{\'e}riodique alors\,}}:
$$N(\phi,P,B) \gg\ll \log B.$$
\end{enumerate}

\begin{rema}\label{remtha}
\begin{itemize}
\item Nous obtenons un r{\'e}sultat plus pr{\'e}cis si nous avons $\deg(\phi)=\deg(\phi^{-1})=d>2$ ainsi que l'une des hypoth{\`e}ses suivantes: 
\begin{itemize}
\item[(H'1)]$\alpha_{arith}(\phi,\pi,V)=\frac{d^2+1}{d},$
\item[(H'2)]$\alpha(\phi,amp)=\frac{d^2+1}{d}$,
\item[(H'3)] $\alpha(\phi,nef)=\frac{d^2+1}{d}.$
\end{itemize}
 alors, nous avons : 
$$N(\phi,P,B)\sim 2\frac{\log B}{\log d}.$$
\item Dans le cas de la dimension 2, nous retrouvons ainsi une  partie des r{\'e}sultats d{\'e}j{\`a} connus pour  le point $(P1)$  la d{\'e}monstration est due {\`a} L. Denis \cite{Denis}, pour le point $(P2)$ la d{\'e}monstration est due {\`a} l'auteure \cite{1art}.
\end{itemize}
\end{rema}



Dans la  partie suivante, nous relions certains indices  au degr{\'e}
dynamique de l'automorphisme, le degr{\'e} dynamique est d{\'e}fini comme suit
(voir par exemple \cite{Sibony}):

\begin{defi}
Soit $\phi$ un automorphisme affine. La limite suivante existe et
d{\'e}fini le degr{\'e} dynamique $\delta(\phi)$:
  $$\delta(\phi):=\inf_{n\geq 1}
deg(\phi^n)^{1/n}$$
\end{defi}

{\`A} partir de r{\'e}sultats de nature arithm{\'e}tique nous obtenons une majoration des indices $\alpha(\phi)$, en fonction du degr{\'e} dynamique de $\phi$ et de $\phi^{-1}$.  

\paragraph*{{\bf{Th{\'e}or{\`e}me B}}}
{\it{Soit $\phi\in \Aut(\mathbb{A}^r(\mathbb{C})).$ Soit $(\pi,V)$ une r{\'e}solution de $\phi$.}}
 {\it{Si $\delta(\phi)=1$, nous posons alors $\delta_1(\phi)=1$,
sinon $\delta_1(\phi)$ est d{\'e}fini comme suit:}}
$\frac{2}{\log \delta_1(\phi)}=\frac{1}{\log
  \delta (\phi)}+\frac{1}{\log \delta (\phi ^{-1})}.$ 
\noindent {\it{Nous avons alors:}}
$$\alpha(\phi, amp)\leq  \delta_1(\phi)+\frac{1}{\delta_1(\phi)},$$
$$\alpha(\phi, nef)\leq  \delta_1(\phi)+\frac{1}{\delta_1(\phi)}.$$
{\it{Nous pouvons d{\'e}finir}} $\alpha (\phi,arith)=\sup_{\pi,V}
\alpha_{arith}(\phi,\pi,V)$ {\it{de plus nous avons:}}
$$\alpha (\phi,arith)\leq  \delta_1(\phi)+\frac{1}{\delta_1(\phi)}.$$

{\it{Nous l'appellerons l'indice effectif arithm{\'e}tique absolu.}}

\begin{rema}
Dans le cas de l'indice effectif arithm{\'e}tique absolu la majoration est optimale dans le sens  o{\`u} nous avons des exemples pour lesquels cette borne est atteinte.
\end{rema}


\noindent
Le th{\'e}or{\`e}me B pour les automorphismes de degr{\'e} dynamique 1; nous indique que l'indice ample est inf{\'e}rieur ou {\'e}gal {\`a} deux; {\`a} partir de r{\'e}sultats de nature arithm{\'e}tique nous montrons que la r{\'e}ciproque de cette assertion est fausse.

La construction d'une r{\'e}solution g{\'e}n{\'e}rique nous permet
l'obtention d'une meilleure borne pour l'indice ample dans tous les
cas en dimension 2 et dans certains cas en dimension sup{\'e}rieure.

\paragraph*{{\bf{Th{\'e}or{\`e}me  C}}}
{\it{Soit $\phi\in \Aut(\mathbb{A}^r)$.}}
{\it{Supposons que $\dim(Z(\phi))=\dim(Z(\phi^{-1}))=0$ alors:}}
$$\alpha(\phi, amp)=0.$$

\begin{rema}
Les conditions de ce th{\'e}or{\`e}me sont toujours v{\'e}rifi{\'e}es en dimension 2. Dans \cite{3art}, nous donnons notamment une autre preuve de ce r{\'e}sultat en dimension 2.
\end{rema}

\noindent Le plan de ce texte est le suivant, nous construisons et {\'e}tudions les invariants g{\'e}om{\'e}triques. Nous utilisons ces indices pour une application arithm{\'e}tique, nous majorons ces invariants, nous montrons que la borne est atteinte pour l'indice effectif, puis nous majorons plus finement l'indice ample.
Dans un souci de lisibilit{\'e} pour le lecteur, nous avons
regroup{\'e}e, dans un appendice, les r{\'e}sultats de dynamique arithm{\'e}tique que nous utilisons.

\paragraph*{Remerciements}   Une partie des r{\'e}sultats {\'e}nonc{\'e}s dans le
th{\'e}or{\`e}me A se trouve dans ma th{\`e}se effectu{\'e}e sous la direction de Marc
Hindry, c'est avec plaisir que je le remercie pour ses conseils et sa
disponibilit{\'e}. Ce travail a pris de l'ampleur lors de mon s{\'e}jour {\`a}
l'universit{\'e} de Ratisbonne (Regensburg), c'est {\'e}galement avec plaisir
que je remercie Uwe Jannsen  et toute l'equipe de th{\'e}orie des nombres
de Ratisbonne pour leur acceuil et cet environnement si
stimulant.C'est {\'e}galement avec plaisir que je remercie Norbert
Schappacher pour ses conseils et avec lui le r{\'e}seau europ{\'e}en
Arithmetic Algebraic Geometry (A.A.G), ainsi que le Max-Planck
Institut f{\"u}r Mathematik, o{\`u} ce travail a {\'e}t{\'e} achev{\'e}.

Enfin, je remercie tr{\`e}s vivement Marco Brunella, pour ses remarques et
questions sur une version pr{\'e}c{\'e}dente de ce texte.


\part{Aspect g{\'e}om{\'e}trique}

\section{Construction des invariants g{\'e}om{\'e}triques}
Dans ce paragraphe, nous montrerons que les indices ample, effectif et nef  sont  des invariants g{\'e}om{\'e}triques des automorphismes de l'espace affine, ce qui revient {\`a} d{\'e}montrer la proposition suivante ainsi que son corollaire.

\begin{prop}\label{eclat}
Soit $\phi \in {\mathbb{A}}^r$ de degr{\'e} $d\geq2$.
Soit $(\pi_1,V_1),(\pi_2,V_2)$ deux r{\'e}solutions de $\phi$.
Nous avons:
 $$\alpha_{max,eff}(\phi,\pi_1,V_1)=\alpha_{max,eff}(\phi,\pi_2,V_2)$$
$$\mbox{et}$$
$$\alpha_{max, nef}(\phi,\pi_1,V_1)=\alpha_{max,nef}(\phi,\pi_2,V_2)$$
\end{prop}


\subsection{Cas de l'indice effectif}
\begin{prop}\label{receclatsup1}
  Soit $\phi$ un automorphisme de ${\mathbb{A}}^r$ de degr{\'e} $d\geq2$.
Soit $\pi:V\rightarrow \mathbb{P}^{r}$ une r{\'e}solution de $\phi$. 
Soit $\pi_1 :V_1\rightarrow V$ un morphisme birationnel. Alors $\pi\circ\pi_1 :V_1\rightarrow \mathbb{P}^{r}$ est encore une r{\'e}solution de $\phi$ et 
$$\alpha_{max,eff}(\pi\circ\pi_1, V_1)=\alpha_{max,eff}(\pi,V).$$
\end{prop}

L'in{\'e}galit{\'e} $\alpha_{max,eff}(\pi,V)\leq \alpha_{max,eff}(\pi\circ\pi_1, V_1)$ est claire, si $\pi_1$ est birationnel r{\'e}gulier il s'agit de prouver l'autre in{\'e}galit{\'e}; de plus il y a {\'e}galit{\'e} si $\pi_1$ est un {\'e}clatement de centre lisse.

La d{\'e}monstration de cette proposition \ref{receclatsup1} se ram{\`e}ne
essentiellement au lemme \ref{eclatsup2} qui nous permet d'obtenir l'autre in{\'e}galit{\'e}.

\begin{lemm}\label{eclatsup2}
Soit $Y$ et $W$ deux vari{\'e}t{\'e}s lisses projectives, avec $W\subset Y$, et $W$
de codimension au moins 2.

\noindent $\underline{Supposons}$~: $\Pi:X\rightarrow Y$ l'{\'e}clatement de $W$
\begin{diagram}
X&\supset &E\\
\dTo^{\Pi}& &\dTo\\
Y&\supset & W\\
\end{diagram}
avec $E\rightarrow W$ fibr{\'e} en $\mathbb{P}^{r}$, et $r+dim(W)=dimX-1$.
On a alors,
$$\Pi^*(\Pic (Y)) \cap \Pic ^+ (X)=\Pi^*(\Pic ^+ (Y)).$$
\end{lemm}

Nous utiliserons le r{\'e}sultat suivant sur les {\'e}clatements.

\paragraph{Rappel}\label{pic} Soit $X$ une vari{\'e}t{\'e} projective lisse, $Y$ une sous-vari{\'e}t{\'e} projective non-singuli{\`e}re de codimension au moins 2. Soit $\pi$ l'{\'e}clatement de $X$ le long de $Y$, $E$ le diviseur exceptionnel associ{\'e}. Nous avons :
$$\Pic(X)= \pi^*(\Pic(Y))\oplus \mathbb{Z}E.$$ 
(voir par exemple \cite{hart} Chapitre II exercice 8.5)

\noindent Les vari{\'e}t{\'e}s sur lesquelles nous travaillons {\'e}tant lisses, nous identifierons les diviseurs de Weil et les diviseurs de Cartier.

\noindent Le symbole $\simeq$ d{\'e}signe l'{\'e}quivalence lin{\'e}aire de diviseurs.

\begin{proof}
Rappelons que $\Pic (X)=\Pi^{\ast} \Pic (Y)\oplus \mathbb{Z}E$ (voir paragraphe \ref{pic}). 
La premi{\`e}re inclusion est assez claire\,: 
$$\Pi^*(\Pic (Y)) \cap \Pic ^+ (X)\supseteq \Pi^*(\Pic ^+ (Y)).$$

\noindent Montrons l'autre inclusion.
Soit $D\in \Pi^*(\Pic (Y)) \cap \Pic ^+ (X)$. 
Comme $D\in \Pic^+ (X)$, on peut {\'e}crire
$$D\simeq \sum m_iD_i +mE$$ avec $D_i$ effectifs irr{\'e}ductibles distincts de  $E$ et $m,m_i\geq 0$.

\noindent Comme $\Pi$ est un isomorphisme hors du support de $E$, pour tout $i$,
il existe $D'_i\in \Pic (Y)$ tel que $\Pi(D_i)=D'_i$. 

\noindent Donc,
$\Pi^*(D'_i)=D_i+e_i E$ et $e_i\in \mathbb{Z}$ pour tout $i$.
\noindent De plus, il existe $\Delta\in \Pic (Y)$ tel que $\Pi^*(\Delta)=D$ avec
$\Delta=\sum n_j \Delta_j$ o{\`u} $\Delta_j$ est irr{\'e}ductible.

\noindent Ainsi, d'une part, nous avons $D=\Pi^*(\Delta)=\sum n_j \Pi^*(\Delta_j)$ et d'autre part
$$D\simeq \sum m_iD_i +mE= \sum m_i \Pi^*(D'_i)+(m-\sum m_ie_i)E$$ 
et par cons{\'e}quent\,:
$$\sum n_j \Pi^*(\Delta_j)\simeq \sum m_i \Pi^*(D'_i)+(m-\sum m_ie_i)E,$$
donc $(m-\sum m_ie_i)=0$;
ou encore $D\simeq \pi^*(\sum m_i D'_i)$,
d'o{\`u} le r{\'e}sultat.
\end{proof}

\begin{proof} (proposition \ref{receclatsup1})
Par hypoth{\`e}se, nous avons\,:
$${(\phi\pi\pi_1)}^\ast (H)+{(\phi ^{-1}\pi\pi_1)}^\ast (H)\simeq \alpha (\pi
 \circ{\pi}_1)^\ast (H) +E_1$$
or $E_1\in  Pic^+ (V_1)$, donc d'apr{\`e}s le lemme
 \ref{eclatsup2}, nous avons $E_1\in \pi_1^*(Pic(V)) \cap Pic^+ (V_1)$
donc, il existe $E$ un diviseur effectif tel que
$E_1=\pi_1^{\ast}(E)$.
Nous avons donc l'autre in{\'e}galit{\'e}.
\end{proof}

\subsection{Cas de l'indice nef}

D'apr{\`e}s \cite{Lazarsfeld} exemple 1.4.3, nous avons:

\begin{lemme}
Soit $f:V\rightarrow W$ un morphisme propre surjectif  de vari{\'e}t{\'e}s lisses projective. Soit $L \in \Pic(W)$, alors:
\begin{center} $L$ est nef si et seulement si $f^*(L)$ est nef. \end{center}
\end{lemme}
Nous en d{\'e}duisons imm{\'e}diatement le lemme suivant :

\begin{lemme}\label{nef}
 Soit $\phi$ un automorphisme de ${\mathbb{A}}^r$ de degr{\'e} $d\geq2$.
Soit $\pi:V\rightarrow \mathbb{P}^{r}$ une r{\'e}solution de $\phi$. 
Soit $\pi_1 :V_1\rightarrow V$ un morphisme birationnel. Alors $\pi\circ\pi_1 :V_1\rightarrow \mathbb{P}^{r}$ est encore une r{\'e}solution de $\phi$ et 
$$\alpha_{max,nef}(\pi\circ\pi_1, V_1)=\alpha_{max,nef}(\pi,V).$$
\end{lemme}

\subsection{Preuve de la proposition \ref{eclat} et de son corollaire}

\begin{proof}
Soit $\pi_1,V_1$,$\pi_2,V_2$ deux r{\'e}solutions de $\phi$. Il existe $W$ une vari{\'e}t{\'e} lisse projective et $f$ un produit fini d'{\'e}clatements tels que: $\psi:=\pi_2^{-1}\pi_1f:W\rightarrow V_2$  est un morphisme de vari{\'e}t{\'e} lisse projective.
Ainsi d'apr{\`e}s la proposition \ref{receclatsup1} dans le cas de l'indice effectif et le lemme \ref{nef} pour l'indice nef nous avons,:
$$\alpha(\pi_1)=\alpha(f\pi_1)=\alpha(\psi\pi_2)\geq\alpha(\pi_2),$$
o{\`u} suivant les cas $\alpha$ d{\'e}signe $\alpha_{max,eff}$ ou bien $\alpha_{max,nef}$.
Pour obtenir l'in{\'e}galit{\'e} inverse, il suffit d'inverser les r{\^o}les des deux r{\'e}solutions.
\end{proof}

La preuve du corollaire s'obtient imm{\'e}diatement {\`a} l'aide du r{\'e}sultat suivant:

\begin{prop}(\cite{Dema} proposition 6.6)
Soit $V$ une vari{\'e}t{\'e} lisse projective. Nous avons:
$$\overline \Pic^a(V)= \Pic^{nef}(V).$$
\end{prop}
 
En effet, nous en d{\'e}duisons le lemme suivant.

\begin{lemm}
  Soit $\phi\in \Aut{\mathbb{A}^r}$ de degr{\'e} au moins deux. Soit $(\pi, V)$ une  r{\'e}solution de $\phi$. Nous avons:
$$\alpha_{max,amp}(\phi, \pi,V)\leq \alpha_{max,nef}(\phi, \pi,V).$$
\end{lemm}

Nous avons donc:
\begin{prop}\label{amp=nef}
Soit $\phi\in \Aut{\mathbb{A}^r}$. Soit $(\pi,V)$ une r{\'e}solution de $\phi$. Nous avons :
$$\alpha_{max, amp}(\phi, \pi,V)\leq\alpha(\phi, nef).$$
\end{prop}

\begin{defi}
Soit $\phi\in \Aut{\mathbb{A}^r}$. L'{\it{indice ample}} existe et est d{\'e}fini comme suit:
$$\alpha(\phi,amp):=\sup\{\alpha_{\max,amp}(\phi,\pi,V)| \pi, V\},$$
\end{defi}
\section{Propri{\'e}t{\'e}s des indices}

\subsection{Stabilit{\'e} par conjugaison lin{\'e}aire}

Nous nous int{\'e}ressons maintenant de mani{\`e}re naturelle {\`a} des questions de stabilit{\'e} des indices.
Nous montrons que les indices ample et effectif sont  invariants par conjugaison par une application lin{\'e}aire. 

\begin{prop}
Soit $\phi,f\in\Aut(\mathbb{A}^{r})$ avec $\deg(f)=1$.

\noindent Nous avons alors,
$$\alpha(\phi, eff)=\alpha( f^{-1}\circ\phi\circ f,eff),$$
$$\alpha(\phi, amp)=\alpha( f^{-1}\circ\phi\circ f,amp).$$
$$\alpha(\phi, nef)=\alpha( f^{-1}\circ\phi\circ f,nef).$$
\end{prop}

\begin{proof}
Soit $\phi\in\Aut(\mathbb{A}^{r})$. Soit $f$ une application lin{\'e}aire de $\mathbb{A}^{r}$. 
Il existe $V$ vari{\'e}t{\'e} projective lisse et $\Pi:V\rightarrow \mathbb{P}^{r}$
morphisme birationnel tels que : $\phi\circ f\circ\Pi$ et $\phi^{-1}\circ f\circ\Pi$  morphismes de  $V$ dans $\mathbb{P}^{r}$.
Comme $f$ est une application lin{\'e}aire, $ f^{-1}\circ\phi\circ f\circ\Pi$ et $ f^{-1}\circ \phi^{-1}\circ f\circ\Pi$ sont des morphismes de $V$ dans $\mathbb{P}^{r}$.
Soit $H$ l'hyperplan {\`a} l'infini.
Le symbole $\simeq$ d{\'e}signe l'{\'e}quivalence lin{\'e}aire de diviseurs.

\begin{eqnarray*}
( f^{-1}\phi f\Pi)^*(H)+ ( f^{-1} \phi^{-1} f\Pi)^*(H)&=& (\phi f\Pi)^*((f^{-1})^*(H))+(\phi^{-1} f\Pi)^*( (f^{-1})^*(H))\\
&\simeq & (\phi f\Pi)^*(H)+(\phi^{-1} f\Pi)^*(H).\\
\end{eqnarray*}
D'o{\`u} le r{\'e}sultat.
\end{proof}

\begin{rema}
La question de l'invariance, par conjugaison par un automorphisme affine, des indices ample, effectif et effectif arithm{\'e}tique  se pose de mani{\`e}re naturelle.
Nous ne sommes pas parvenue {\`a} trancher. N{\'e}anmoins le th{\'e}or{\`e}me B, nous
permet de d{\'e}finir l'{\it{indice ample conjugu{\'e}}}, l'{\it{indice nef conjugu{\'e}}} et l'{\it{indice
    effectif arithm{\'e}tique absolu conjugu{\'e}}} :
 $$\alpha_{conj,}(\phi, amp)=\sup_{f\in \Aut (\mathbb{A}^r)}
 \alpha(f^{-1}\phi f, amp),$$
$$\alpha_{conj,}(\phi, nef)=\sup_{f\in \Aut (\mathbb{A}^r)}
 \alpha(f^{-1}\phi f, nef),$$
$$\alpha_{conj,}(\phi, arith )=\sup_{f\in \Aut (\mathbb{A}^r)}
 \alpha(f^{-1}\phi f, arith ).$$
En effet, le degr{\'e} dynamique est invariant par conjugaison par un automorphisme de l'espace affine. 
\end{rema}

\subsection{Relation entre les indices}

Nous avons une in{\'e}galit{\'e} dans le cas g{\'e}n{\'e}ral. 
\begin{prop}
Soit $\phi\in \Aut(\mathbb{A}^r)$. Alors :
$$\alpha(\phi,amp)\leq \alpha(\phi,eff).$$
\end{prop}
 
Ceci d{\'e}coule imm{\'e}diatement du fait que $\Pic^{a}_{\mathbb{R}}(V)\subset \Pic^+_{\mathbb{R}}(V)$, voir par exemple \cite{hs} F.5.4.


\begin{rema}
Soit $\phi\in \Aut(\mathbb{A}^r)$. Par d{\'e}finition nous avons:
$$\alpha (\phi, arith)\leq \alpha(\phi,eff).$$
Dans tous les exemples que nous connaissons, nous avons {\'e}galit{\'e}.  
\end{rema}

\begin{prop}
Soit $\phi\in \Aut(\mathbb{A}^{r})$, $r\geq 2$ et  $\deg(\phi)=1$.
Alors,
$$\alpha(\phi, amp)=\alpha(\phi, nef) =\alpha(\phi, eff)=2.$$
\end{prop}

\begin{proof}
De l'hypoth{\`e}se $\deg(\phi)=1$, nous d{\'e}duisons imm{\'e}diatement que $\phi$ et $\phi^{-1}$ s'{\'e}tendent en des morphismes de $\mathbb{P}^{r}$, de l{\`a} nous avons $\alpha(\phi, eff)=2.$
De plus, $\Pic(\mathbb{P}^{r})=\mathbb{Z}$, or d'apr{\`e}s \cite{Lazarsfeld} ex. 1.23, tout diviseur effectif non nul est ample, d'o{\`u} le r{\'e}sultat. 
\end{proof}

\subsection{G{\'e}n{\'e}ralisation}
Il est possible de consid{\'e}rer les indices ample et effectif comme {\'e}tant les invariants g{\'e}m{\'e}triques de l'application birationnelle $\tilde \phi:\mathbb{P}^r\rightarrow \mathbb{P}^r$ associ{\'e}e {\`a} $\phi$. Ainsi en vient-on de mani{\`e}re naturelle {\`a} d{\'e}finir les {\it{indices nef et effectif}} pour toute application birationnelle d'une vari{\'e}t{\'e} lisse projective $V$ de dimension au moins 2
munie d'une polarisation $L$, o{\`u} $L$ jouera le r{\^o}le de l'hyperplan {\`a} l'infini $H$.
Nous noterons alors ces indices:
$\alpha(\phi,eff,L)$ et $\alpha(\phi,nef,L)$.
La preuve de la proposition \ref{receclatsup1} s'adapte {\it{mutatis mutandis} }{\`a} ce cas de figure.


\part{Dynamique arithm{\'e}tique}

Le th{\'e}or{\`e}me A s'obtient {\`a} l'aide des propri{\'e}t{\'e}s fonctorielles des hauteurs et de lemmes combinatoires.

\section{Propri{\'e}t{\'e}s fonctorielles des  hauteurs}

La d{\'e}monstration de la proposition suivante repose essentiellement sur les propri{\'e}t{\'e}s fonctorielles des hauteurs voir par exemple \cite{hs} B.3.

\begin{prop}\label{weil} Soit $\phi\in \Aut(\mathbb{A}^r)$.
Alors, il existe
  $c\in \mathbb{R}^+$ tel que pour tout point Q de ${\mathbb{A}}^r({\overline{\mathbb{Q}}})$, on a:
  $$h(\phi(Q))+h(\phi^{-1}(Q))\geq \alpha h(Q)-c \eqno{(1)},$$
o{\`u} $\alpha$ d{\'e}signe suivant les cas :
\begin{itemize}
\item $\alpha\in\mathbb{Q}$ tel que $2<\alpha\leq \alpha_{arith}
  (\phi, \pi,V)$ avec $(\pi,V)$ une r\'esolution de $\phi$,
\item $\alpha\in\mathbb{R}$ tel que $2<\alpha<\alpha(\phi,amp)$,
\item $\alpha\in\mathbb{R}$ tel que $2<\alpha<\alpha(\phi,nef)$.
\end{itemize}
\end{prop}

\begin{rema}
La constante $c$ d{\'e}pend la r{\'e}solution choisie. 
\end{rema}

\begin{proof} Soit $(\pi,V)$ une r{\'e}solution de $\phi$.
Nous noterons $\psi:=\phi\pi$ et $\psi':={\phi}^{-1}\pi$.
Si D est un diviseur de Weil sur V (vari{\'e}t{\'e} projective non singuli{\`e}re), on notera
 $h_{V,D}$ la hauteur sur V associ{\'e}e au diviseur D.
\begin{itemize}
\item Nous consid{\'e}rons tout d'abord le cas de l'indice effectif arithm{\'e}tique.

\noindent Par hypoth{\`e}se, il existe une r{\'e}solution $(\pi,V)$ et $E\in \Cplus(V)$ pour laquelle:
  $$(\phi\pi)^*(H)+({\phi}^{-1}\pi)^*(H)- \alpha \simeq E\ \eqno{(2)},$$ 
o{\`u} $\simeq$ d{\'e}signe l'{\'e}quivalence lin{\'e}aire de diviseurs, et il existe $m>0$ tel que pour tout $(Q\notin H)$ nous avons $\pi^{-1}(Q)\notin \B_m(E) $.

$$h_{V,\psi^*(H)}+h_{V,{\psi'}^*(H)}=h_{V,\psi^*(H)+{\psi'}^*(H)}+O(1)$$

$$h_{V,\psi^*(H)}+h_{V,{\psi'}^*(H)}=h_{V,\alpha\pi^*(H)+E}+O(1)$$

$$h_{V,\psi^*(H)}+h_{V,{\psi'}^*(H)}=\alpha h_{V,\pi^*(H)}+h_{V,E}+O(1)$$

Pour tout $Q\in {\Bbb A}^r,(Q\notin H)$, nous  noterons $\tilde Q =\pi^{-1}(Q)$
par hypoth{\`e}se il existe $m>0$ tel que $\tilde Q \notin \B_m(E) $.

Ainsi, pour
    tout $Q\in \mathbb{A}^r(\bar{Q}),$
$$h_{V,\psi^*(H)}(\tilde Q)+h_{V,{\psi'}^{*}(H)}(\tilde Q)=\alpha
h_{V,\pi^{*}(H)}(\tilde Q)+h_{V,E}(\tilde Q)+O(1).$$

Et par cons{\'e}quent,
$$h(\psi(\tilde Q))+h(\psi^{-1}(\tilde Q))=\alpha h(\pi(\tilde Q ))+h_{V,E}(\tilde Q)+O(1)$$
or,

$Q=\pi(\tilde Q)$, donc $\psi(\tilde Q)=\psi\circ\pi(Q)=\phi(Q)$
d'o{\`u}
$$h(\phi(Q))+h(\phi^{-1}(Q))=\alpha h(Q)+h_{V,E}(\tilde Q)+O(1)$$
de plus, E est un diviseur effectif, $m>0$ et $\tilde Q\notin \B_m(E)$,
 donc il existe $c>0$ tel que:
$$h(\phi(Q))+h(\phi^{-1}(Q))\geq\alpha h(Q)-c.$$

\item Cas des indices amples et nef.
Si $E$ est un diviseur ample, alors il existe $c$ tel que pour tout $P$,  $h_E(P)\geq -c$, nous pouvons donc adapter la preuve pr{\'e}c{\'e}dente.
Le cone ample est ouvert d'o{\`u} la n{\'e}c{\'e}ssit{\'e} de consid{\'e}rer $\alpha<\alpha(\phi, amp)$.
\end{itemize}

\end{proof}

\section{Techniques combinatoires}

L'{\'e}nonc{\'e} suivant figure, sans d{\'e}monstration dans \cite{moi}.

\begin{prop}\label{tech1}
Soit $\phi$ un automorphisme de ${\mathbb{A}}^r$ de degr{\'e} $d\geq 2$
  poss{\'e}dant un inverse de degr{\'e} $d'$.

\noindent  Supposons qu'il existe
  $c\in {\mathbb{R}}^+$, $\alpha>2$ tel que pour tout point Q de ${\mathbb{A}}^r(\overline{{\mathbb{Q}}})$\,:
  $$h(\phi (Q))+h(\phi^{-1}(Q))\geq \alpha h(Q)-c. \eqno(1)$$

\begin{enumerate}
\item Notons $Per(\phi)$ l'ensemble des points p{\'e}riodiques de
  $\phi$. Nous avons alors\,:
$$Per(\phi)\subset \Big\{ Q\in {\mathbb{A}}^r(\overline{\mathbb{Q}}) \mid h(Q)\leq  \frac{c}{\alpha-2}\Big\}.$$
\item Soit $a(\alpha)=\frac {\alpha +({\alpha
      }^2-4)^{\frac{1}{2}}}{2}$. Pour tout point $P$ non-p{\'e}riodique, on a\:
$$\Big( \frac{1}{\log d}+\frac{1}{\log d'}\Big)\leq \liminf
\frac{N(\phi,P,B)}{\log B}\leq  \limsup \frac{N(\phi,P,B)}{\log B}\leq \frac{2}{\log (a(\alpha))}. \eqno(2) $$
\end{enumerate}
\end{prop}

\begin{rema}
La minoration dans la relation $(2)$ est vraie pour tout automorphisme affine $\phi$.

\noindent De plus, si on note $\delta$ (resp.$\delta '$) le degr{\'e}
dynamique de $\phi$ (resp. de $\phi ^{-1}$); on peut remplacer, dans
cette m{\^e}me minoration,  $d$ et $d'$ respectivement par $\delta$  et $\delta '$ (voir \cite{1art} proposition 2.8 ).  
\end{rema}

La d{\'e}monstration de la proposition  \ref{tech1} est essentiellement technique et utilise les
propri{\'e}t{\'e}s de bases des hauteurs ainsi que des techniques classiques
pour les suites. 
La d{\'e}monstration de la proposition \ref{tech1} est largement inspir{\'e}e des articles de J. Silverman \cite{Sil1} et L. Denis \cite{Denis}.

\subsection{Lemmes techniques}

\begin{lemm}\label{lem1}
Soit $\phi$ un automorphisme de ${\mathbb{A}}^r$, v{\'e}rifiant la
condition $(1)$ alors, pour tout point $P$ de  ${\mathbb{A}}^r$ nous avons\,:

 $$h({\phi}^n
 (P))+h({\phi}^{-n}(P))\geq(a-a^{-1})(a^n+a^{-n})\Big(h(P)-\frac{c}{\alpha
   -2}\Big),$$
avec  $a=\frac {\alpha+(\alpha^2-4)^{\frac{1}{2}}}{2}$.
\end{lemm}

L'objectif est d'obtenir pour tout point $P$ de l'espace affine une minoration de $h({\phi}^n
 (P))+h({\phi}^{-n}(P))$ en fonction de $h(P)$ et de constantes. Pour ce
 faire, nous utilisons une suite r{\'e}elle associ{\'e}e de mani{\`e}re
 naturelle {\`a} l'in{\'e}galit{\'e} $(1)$, {\`a} savoir une suite qui v{\'e}rifie le
 cas d'{\'e}galit{\'e} avec des conditions initiales qui nous conviennent.
Ensuite {\`a} l'aide d'un argument de sommation t{\'e}l{\'e}scopique,
 puis de minorations {\'e}l{\'e}mentaires, nous obtenons le r{\'e}sultat escompt{\'e}.

Pour la preuve de ce lemme nous renvoyons le lecteur {\`a} \cite{1art} lemme 3.3.





Ce lemme nous permet d'obtenir la premi{\`e}re partie de la proposition \ref{tech1}. Pour avoir des informations sur l'orbite des points non-p{\'e}riodiques, nous utilisons le lemme suivant\,:

\begin{lemm}\label{lem2}
  Soit $\phi$ un automorphisme du plan affine v{\'e}rifiant
  la relation $(1)$. Si $P$ n'est pas un point p{\'e}riodique, alors les suites
  $(h(\phi^n(P)))_{n\in {\Bbb N}}$ et $(h(\phi^{-n}(P)))_{n\in {\Bbb
      N}}$ sont croissantes {\`a} partir d'un certain rang.
\end{lemm}

\begin{proof}
   Supposons qu'il existe $m$ tel que
  $h(\phi^m(P))\leq h(\phi^{m+1}(P))$ et  $h(\phi^{m+1}(P))\geq
  h(\phi^{m+2}(P))$. On a alors,
  $$2h(\phi^{m+1}(P))\geq h(\phi^m(P))+h(\phi^{m+2}(P))$$ d'o{\`u} en
  utilisant la propri{\'e}t{\'e} (1) avec $Q=\phi^{m+1}(P)$
  $$2h(\phi^{m+1}(P))\geq \alpha h(\phi^{m+1}(P))-c$$ et par
  cons{\'e}quent,
  $$c\geq (\alpha -2)h(\phi^{m+1}(P)).$$

Cependant, si $P$ est non-p{\'e}riodique $\#\{m\in\mathbb{Z} \mid
  h(\phi^{m+1}(P))\leq \frac{c}{\alpha -2} \}$ est fini.
D'o{\`u} la contradiction pour $m$ assez grand.
On proc{\`e}de de la m{\^e}me mani{\`e}re pour d{\'e}montrer la croissance {\`a} partir
  d'un certain rang de la suite de terme g{\'e}n{\'e}ral $(h(\phi^{-n}(P)))_{n\in {\mathbb N}}$.
\end{proof}

\begin{rema}
{\`A} partir d'un certain rang la suite $(H_n)_{n\in {\Bbb N}} $ est
  {\'e}galement 
  croissante.
\end{rema}

Le lemme suivant nous donne
  une minoration (ind{\'e}pendante du point consid{\'e}r{\'e}), de la limite inf{\'e}rieure de la hauteur de it{\'e}r{\'e}s.

\begin{lemm}\label{lem3}
Soit $\phi$ un automorphisme de l'espace affine $\mathbb{A}^r$
v{\'e}rifiant la condition $(1)$. Notons $a=\frac {\alpha+\sqrt{\alpha^2-4}}{2}$.

\noindent Alors pour tout point $P$ non-p{\'e}riodique nous avons:
$$\liminf_{n\rightarrow +\infty}\frac{\log h({\phi}^n(P))}{n\log
  a}\geq 1\ .$$

\end{lemm}

Sa d{\'e}monstration, que voici, n{\'e}cessite l'utilisation des lemmes \ref{lem1} et
\ref{lem2}.

\begin{proof}
Comme $P$ n'est pas un point p{\'e}riodique,
    il existe $n_0(P)$, tel que $h(\phi^s(P))$ est une suite croisssante ({\it cf.}
    lemme \ref{lem2}) pour tout $s>n_0(P)$, et tel que, $h(\phi^{s}P)\geq
    \frac{ac}{2(a-1)}$ de l{\`a},

    $h(\phi^{s+k} (P))\geq h(\phi^{s-k} (P))$, pour $s>n_0+k$ avec $k>0$

    d'o{\`u} pour $s> n_0+k$
$$
   2h(\phi^{s+k} P)\geq h(\phi^{s-k} P)+h(\phi^{s+k} P)= h(\phi^{-k}(\phi^{s}
    P))+h(\phi^{k}(\phi^{s} (P) )
$$

    donc, d'apr{\`e}s le lemme \ref{lem1}
\begin{xxalignat}{2}
    2h(\phi^{s+k} (P))&\geq  c+(a-a^{-1})(a^k+a^{-k})(h(\phi^{s}
    P)-\frac{c}{2})& \qquad {\mbox{(lemme \ref{lem1})}}\\
   &\geq (a-a^{-1})(a^k+a^{-k})(h(\phi^{s}(P))-\frac{c}{2}) &\qquad {\mbox{($c>0$ )}} \\
& \geq (a-a^{-1})(a^k+a^{-k})\frac{1}{a}h(\phi^{s}(P)) &\qquad (car\  h(\phi^{s}P)\geq \frac{ac}{2(a-1)}).
\end{xxalignat}
Posons, $c_1=\frac{1}{2}(1-\frac{1}{a^2})$. Nous avons alors,
$$h(\phi^{s+k} (P))\geq c_1 a^k h(\phi^{s}P);$$
    d'o{\`u} en appliquant $m$ fois cette in{\'e}galit{\'e}\,:
    $$h(\phi^{s+mk} (P)\geq c_1^m a^{km}h(\phi^{s}(P)\ pour\ tout\ 
    s>n_{0}+k\ et\ tout\ m\geq 0.$$ 
De l{\`a}, 
\begin{xxalignat}{2}
   \frac{\log h(\phi^{s+mk} (P))}{s+mk}&\geq (1-\frac{s}{s+mk})\log a+ \frac{m\log c_1 +\log h(\phi^{s}(P)}{s+mk},
\end{xxalignat}
d'o{\`u} le r{\'e}sultat par passage {\`a} la limite.

\end{proof}

\subsection{Preuve de la proposition \ref{tech1}}

\begin{proof} (de la proposition \ref{tech1})
Le premier point r{\'e}sulte d'une simple application du lemme \ref{lem1}.
\begin{enumerate}
\item D'apr{\`e}s le lemme \ref{lem1}   
    $$h(\phi^n(P))+h(\phi^{-n}(P))\geq
    (a-a^{-1})(a^n+a^{-n})(h(P)-\frac{c}{\alpha-2})+\frac{2c}{\alpha-2}$$

    Soit P un point p{\'e}riodique. Il existe $n\in \mathbb{N}$ tel que $\phi^n(P)=P$, ainsi
    $$\frac{c}{\alpha-2}[(a-a^{-1})(a^n+a^{-n})-2]\geq h(P)
    [(a-a^{-1})(a^n+a^{-n})-2]$$
d'o{\`u} le r{\'e}sultat.

\vspace{1cm}
Pour obtenir un encadrement de la taille des orbites, nous proc{\'e}dons
tout d'abord {\`a} un encadrement de  $\log(h(\phi^n(P)))$, puis
comptons les points de hauteur born{\'e}e de l'orbite.

\item Taille de l'orbite

Pour tout $P\in A^r({\bar{\mathbb{Q}}})$, $h(\phi^nP)\leq
    d^n(h(P)+O(1))$ en effet $\phi$ est une application
    rationnelle de degr{\'e} $d$. Comme $P$ n'est pas un point p{\'e}riodique, nous avons,

    $$\limsup_{n\rightarrow\infty }\frac{\log h(\phi^nP)}{n\log d}\leq
    \limsup_{n\rightarrow\infty }\frac{n\log d+\log (h(P))}{n\log
      d}=1;$$
    par cons{\'e}quent, pour $n>n_0$ nous avons \,:
   $$ n \log a \leq \liminf_{n\rightarrow\infty } \log h(\phi^nP)\leq \limsup_{n\rightarrow\infty } \log h(\phi^nP)\leq n\log d.$$
    Le r{\'e}sultat pour $n\rightarrow -\infty$ s'obtient de mani{\`e}re analogue car par
    hypoth{\`e}se $\phi^{-1}$ est un automorphisme affine dont
    l'application rationnelle associ{\'e}e est de degr{\'e} $d'$, ainsi\,
    pour tout $\varepsilon>0$ il existe $n_o$ tel que pour tout n
    v{\'e}rifiant $\mid \!n\!\mid >n_0$\,:
$$(1-\varepsilon)\mid\! n\!\mid \log a \leq \log h(\phi^nP)\leq (1+\varepsilon)\mid
    \!n\!\mid \log (\max(d,d')) $$ pour tout $\mid n\mid\geq n_0$. 
 Nous avons,
    $$N(\phi,P,B)=\#\{n\in \mathbb{Z}\quad  \mid
    \quad \log(h(\phi^nP))\leq \log B\}$$ on a:
    $$\frac{ \log B}{(1+\varepsilon)(\frac{1}{\log d} +\frac{1}{\log d'})}-2n_0-2\leq N(\phi,P,B)\leq \frac{2 \log
      B}{(1-\varepsilon)\log a}+1$$ 
 d'o{\`u} le r{\'e}sultat.

\end{enumerate}
\end{proof}

\section{D{\'e}monstration du th{\'e}or{\`e}me A}

La d{\'e}monstration du th{\'e}or{\`e}me A et de la remarque \ref{remtha} se d{\'e}duit imm{\'e}diatement des propositions \ref{weil} et \ref{tech1}.
Nous traitons le cas de l'indice ample sachant que la preuve pour l'indice effectif arithm{\'e}tique est essentiellement la m{\^e}me.

Le cas d'{\'e}quivalence \ref{remtha} s'obtient comme suit\,:

\begin{proof}[\bf{D{\'e}monstration du th{\'e}or{\`e}me A}]
 Par hypoth{\`e}se $\alpha (\phi,amp)=\frac{d^2+1}{d}$, de plus il existe $n_1\in \mathbb{N}\setminus\{0\}$ tel que pour tout $n>n_1$ nous avons $\alpha (\phi,amp)- \frac{1}{n}> 2$
Or $a_n=\frac
    {\alpha(\phi,amp)-\frac{1}{n} +\sqrt{{(\alpha(\phi,amp)}-\frac{1}{n}) ^2-4}}{2}$, donc $\lim_{n\rightarrow \infty} a_n=d$ de plus $d=d'$, do{\`u}\,:
$$\lim_{\mid\!n\!\mid \rightarrow\infty }\frac{\log h(\phi^nP)}{n\log d}=1$$
 Ainsi, pour tout $\varepsilon$,il existe $n_0(P,\varepsilon)$ tel que:

   \noindent $(1-\varepsilon)\mid n\mid\leq \log_dh(\phi^nP)\leq (1+\varepsilon)\mid
    n\mid$ pour tout $\mid n\mid\geq n_0$.  Nous avons alors,
   
    $$\frac{2 \log_d B}{1+\varepsilon}-2n_0-2\leq N(\phi,P,B)\leq \frac{2 \log_d
      B}{1-\varepsilon}+1$$ d'o{\`u}
    $$\frac{1}{1+\varepsilon}\leq \liminf_{B\rightarrow\infty}\frac{N(\phi,P,B)}{2
      \log_d B}\leq \limsup_{B\rightarrow\infty}\frac{N(\phi,P,B)}{2 \log_d
      B}\leq \frac{1}{1-\varepsilon}$$ d'o{\`u} le r{\'e}sultat.
\end{proof}

\section{Illustration du th{\'e}or{\`e}me A} 

Nous illustrons le th{\'e}or{\`e}me  A gr{\^a}ce aux calculs des invariants g{\'e}om{\'e}triques. 

\subsection{Valeurs num{\'e}riques des indices effectif et effectif arithm{\'e}tique}

Dans ce paragraphe, nous donnons les valeurs des deux invariants
g{\'e}om{\'e}triques pour diff{\'e}rents exemples, dans tous ces exemples $a$ est
une constante non nulle. J. Silverman a d{\'e}termin{\'e} la valeur de
l'indice effectif pour l'un de ces exemples, l'application de H{\'e}non
g{\'e}n{\'e}ralis{\'e}e \cite{Sil1}. Pour calculer ces indices nous avons
utilis{\'e} la r{\'e}solution canonique \cite{3art}. 
Le d{\'e}tail des calculs est omis. Dans les exemples ci-dessous $a$ est toujours une constante non nulle.





\vspace{1cm}

\noindent\begin{tabular}{|c||c|c|c|c|c|}
\hline
Application & degr{\'e} & degr{\'e}  & indice& indice& indice  \\
 &  alg{\'e}brique &  dynamique & effectif& arithm{\'e}tique & ample \\
\hline
$\phi(x,y)=$& $\deg(\phi)$ & $\delta(\phi)$ & $\alpha(\phi,eff)$&  $\alpha_{arith,can}(\phi)$&  $\alpha(\phi,amp)$\\
\hline
\hline
$(y,y^2+b+ax)$&2   & 2& $ \frac{5}{2}=2+\frac{1}{2}$ \cite{Sil1}& $ \frac{5}{2}$& $\leq 0$\\
\hline
$(y+ax^3,x)$ &3 & 3 & $\frac{10}{3}=3+\frac{1}{3}$& $\frac{10}{3}$& $\leq 0$\\
\hline
$(y+ax^4,x)$&4 &4& $\frac{17}{4}=4+\frac{1}{4}$&$\frac{17}{4}$ & $\leq 0$\\
\hline 
\end{tabular}

\vspace{1cm}

\begin{rema} Pour les exemples de degr{\'e} dynamique diff{\'e}rents de 1, nous obtenons la borne maximale du th{\'e}or{\`e}me C.
Sur les calculs effectu{\'e}s, nous observons :
$\alpha(\phi,eff)=\alpha_{arith,conj}(\phi,\pi,V)=\delta(\phi)+\frac{1}{\deg(\phi)}.$
Une question naturelle se pose, les indices amples et effectifs ne seraient-ils pas rationnels?
\end{rema}

\subsection{Taille des orbites}

 Dans les exemples ci-dessous $a$ est toujours une constante non nulle.
En appliquant ce th{\'e}or{\`e}me A nous obtenons que l'ensemble des points p{\'e}riodiques des applications suivantes est de hauteur born{\'e}e sur $\overline{\mathbb{Q}}$, et pour les tailles d'orbites des points non-p{\'e}riodiques, nous avons\,:

\vspace{1cm}

\hspace{1cm}\begin{tabular}{|c||c|c|}
\hline
Application & degr{\'e} dynamique & taille de l'orbite \\
$\phi(x,y)=$ & $\delta(\phi)$ & $N(\phi,P,B)$\\
\hline
$(y,y^2+b+ax)$   & 2& $\sim \frac{2\log B}{\log 2}$ \cite{Sil1}\\
\hline
$(y+ax^3,x)$  & 3 & $\sim\frac{2\log B}{\log 3}$\\
\hline
$(y+ax^4,x)$ &4& $\sim\frac{2\log B}{\log 4}$\\
\hline
\end{tabular}

\vspace{1cm}

Un autre r{\'e}sultat que nous avons obtenu \cite{1art} concernant les automorphismes r{\'e}guliers (voir Annexe A th{\'e}or{\`e}me \ref{reg}), nous permet de retrouver ces r{\'e}sultats.

\section{Cas o{\`u} le th{\'e}or{\`e}me A n'est pas optimal}

Dans certains cas, le th{\'e}or{\`e}me A ne donne pas un r{\'e}sultat optimal.

Nous consid{\'e}rons le cas de l'application de Nagata-``tordue'' qui est d{\'e}finie comme suit:
\[
\begin{array}{lccl}
\phi : & \mathbb{A}^3 &\to& \mathbb{A}^3\\
& (X,Y,Z) &\mapsto &
\left\{
\begin{array}{l}
Y-2(YZ+X^2)X-(YZ+X^2)^2Z\\
X+(YZ+X^2)Z\\
Z,\\
\end{array}
\right. \\
\end{array}     
\]
 Le th{\'e}or{\'e}me \ref{nagt} et la contrapos{\'e}e du th{\'e}or{\`e}me
A nous permettent de conclure que si l'indice de
l'application de Nagata-``tordue'' v{\'e}rifie  respectivement $\alpha (\phi, nef)>2$ et $\alpha (\phi,arith>2$ , alors les deux conditions suivantes sont v{\'e}rifi{\'e}es:
\begin{itemize}
\item l'indice ample v{\'e}rifie  $\alpha (\phi, nef) < \frac{26}{5}$,
\item soit l'indice effectif arithm{\'e}tique v{\'e}rifie  $\alpha  (\phi, arith)< \frac{26}{5}$.
\end{itemize}

Cet exemple nous montre que dans ce cas le th{\'e}or{\`e}me A ne permet pas, dans ce cas, d'obtenir un r{\'e}sultat optimal. Il ne donnerait au mieux que ``comparable {\`a}'' $\log B$, alors que le th{\'e}or{\'e}me \ref{nagt} nous d{\'e}crit pr{\'e}cisement la taille des orbites.


\part{Majoration des indices}

\section{Preuve et illustration du th{\'e}or{\`e}me B}

Pour d{\'e}montrer le th{\'e}or{\`e}me B, nous consid{\'e}rons deux cas suivant que le degr{\'e} dynamique est ou non {\'e}gal {\`a} 1. 
Dans les deux cas, nous utilisons des r{\'e}sultats de dynamique arithm{\'e}tique puis des arguments sp{\'e}cialisation pour le passage de $\overline{\mathbb{Q}}$ {\`a} $\mathbb{C}$.

\subsection{Automorphismes de degr{\'e} dynamique 1}

Nous utilisons des r{\'e}sultats de dynamique arithm{\'e}tique que nous avons obtenus pr{\'e}c{\'e}demment \cite{1art}.

\noindent Le th{\'e}or{\`e}me B pour les automorphismes affines de degr{\'e} dynamique 1 s'obtient {\`a} partir de la proposition 
\ref{hdyn}, de la contrapos{\'e}e  du th{\'e}or{\`e}me A et d'arguments
de sp{\'e}cialisation pour le passage de $\mathbb{\overline{Q}}$ {\`a} $\mathbb{C}$.

\begin{rema}
La majoration est optimale, il suffit en effet de consid{\'e}rer l'identit{\'e} qui est d'indice 2. 
\end{rema}

\noindent La situation produit (voir Annexe paragraphe A.3) nous permet de montrer que la r{\'e}ciproque de l'assertion du  th{\'e}or{\`e}me B pour les automorphismes affines de degr{\'e} dynamique 1 est fausse.

Nous pouvons par exemple consid{\'e}rer l'exemple suivant.
Soit $\phi_1\in \Aut(\mathbb{A}^{2})$ une application triangulaire. 
Soit $\phi_2$  une application de H{\'e}non g{\'e}n{\'e}ralis{\'e}e, c'est-{\`a}-dire $\phi_2(x,y)=(p(x)-ay,x)$ avec $p$ application polynomiale de degr{\'e} au moins 2 et $a$ une constante non nulle.

\noindent L'automorphisme $\Phi$ a pour
expression\,: 
$$\Phi(x,y,z,t)=(\alpha x+q(y),\beta y+\gamma,p(z)-at,z),$$
avec $a\alpha\beta\neq 0$ et $p$ application polynomiale de degr{\'e} au moins 2. 

Des r{\'e}sultats de dynamique arithm{\'e}tique pour la situation produit voir Annexe A paragraphe A.3, nous d{\'e}duisons\,:
\noindent si $P=(x,y,z,t)$ avec $P_1=(x,y)\notin Per(\phi_1,k)$
  et $P_2=(z,t)\in  Per(\phi_2,k)$ alors,
 $$N(\Phi,P,B)=N(\phi_1,P_1,B).$$ 
Suivant le type d'automorphisme
 triangulaire et le point consid{\'e}r{\'e}s, diff{\'e}rents cas se pr{\'e}sentent. 
Les tailles des orbites sont donc des
 fonctions qui asymptotiquement se comportent comme $cB$ ou $ce ^B$.
Ainsi, de la contrapos{\'e}e du th{\'e}or{\`e}me A, nous d{\'e}duisons $\alpha(\Phi,nef)\leq 2$ de m{\^e}me que $\alpha(\Phi, arith)\leq 2$.
De plus, nous avons $\delta(\Phi)=\delta(\phi_2)>1$, en effet une
application de H{\'e}non g{\'e}n{\'e}ralis{\'e}e $\phi$ v{\'e}rifie
$\deg(\phi)=\delta(\phi)$ (voir par exemple \cite{Sibony}).
\noindent La r{\'e}ciproque du th{\'e}or{\`e}me C (pour les automorphismes de
degr{\'e} dynamique 1) est donc fausse.

\subsection{Automorphismes de degr{\'e} dynamique  diff{\'e}rent de 1}
Soit $\phi$ un automorphisme de l'espace affine de degr{\'e} dynamique $\delta(\phi)>1$.
Nous devons distinguer deux cas suivant que l'indice nef ou l'indice effectif arithm{\'e}tique est ou non inf{\'e}rieur ou {\'e}gal {\`a} 2. Si $\alpha(\phi, nef)\leq 2$ ou $\alpha (\phi, arith)\leq 2$, l'assertion est clairement v{\'e}rifi{\'e}e.

Dans le cas contraire, si $\alpha(\phi,nef)> 2$ ou $\alpha (\phi,arith)> 2$ le th{\'e}or{\`e}me B s'obtient essentiellement {\`a} partir de l'in{\'e}galit{\'e} 2 de la
proposition \ref{tech1} modifi{\'e}e par la remarque 3.3.
La borne est optimale, elle est en effet atteinte, 
dans certains calculs explicites de l'indice que nous avons
effectu{\'e}s (voir paragraphe 16). En effet sur ces exemples, nous avons
{\'e}galit{\'e} des degr{\'e}s dynamiques et alg{\'e}briques, et de plus
l'application et son inverse ont m{\^e}me degr{\'e}, la borne
sup{\'e}rieure devient alors $d+\frac{1}{d}.$  

\begin{rema}
Il est  possible de d{\'e}finir de mani{\`e}re analogue une borne en fonction des degr{\'e}s alg{\'e}briques et nous avons alors $\delta_1(\phi)\leq \deg_1 (\phi)$.
\noindent Le nombre $\delta_1(\phi)$ est clairement invariant par conjugaison par un automorphisme affine.
\end{rema}

Il serait int{\'e}ressant d'avoir une d{\'e}monstration directe du th{\'e}or{\`e}me B.

\subsection{Illustrations du th{\'e}or{\`e}me B}

Dans ce paragraphe, nous donnons diff{\'e}rents exemples d'automorphismes de degr{\'e} dynamique 1, qui nous permettent donc d'illustrer par des exemples le th{\'e}or{\`e}me C pour les automorphismes de degr{\'e} dynamique 1.

Les automorphismes de degr{\'e} dynamique 1 sont nombreux. Il s'agit notamment des applications triangulaires et {\'e}l{\'e}mentaires.

\begin{defi}
Soit $\phi \in\Aut(\mathbb{A}^r)$, v{\'e}rifiant
\begin{eqnarray*}
\phi(X_{1},\dots,X_{r})&=&(a_1X_{1}+F_1(X_{2},\dots,X_{r}),a_2X_2+F_2(X_{3},\dots,X_{r}),\\
& &\quad \dots,a_{r-1}X_{r-1}+F_{r-1}(X_{r}),a_rX_r+F_r);
\end{eqnarray*}
avec   $F_i\in k[X_{i+1},\dots,X_{r}]$ pour $1\leq i\leq r-1$  et  $a_i\in
k^*$;
$\phi$ est dite {\it{triangulaire}}.
\end{defi}
Les applications triangulaires jouent un r{\^o}le important dans la classification des automorphismes du plan affine (voir \cite{Jung},\cite{Kulk} et \cite{FM}).
\begin{defi}
Soit $\phi \in\Aut(\mathbb{A}^r)$, v{\'e}rifiant
\[
\begin{array}{llcl}
\phi: & \mathbb{A}^r & \to & \mathbb{A}^r \\
 &(X_1,..,X_i,...,X_r) & \mapsto &
    (X_1,..,X_{i-1},X_i+P(X_1,..,\hat{X_i},.,X_r),X_{i+1},...X_r)
\end{array}     
\]
 avec $P$ polyn{\^o}me; $\phi$ est dite {\it{{\'e}l{\'e}mentaire}}.
\end{defi}
Les applications {\'e}l{\'e}mentaires figurent {\'e}galement dans les travaux sur les automorphismes affines \cite{Essen}.
L'application d'Anick est {\'e}galement de degr{\'e} dynamique 1, elle est d{\'e}finie comme suit:
\[
\begin{array}{lccl}
\sigma : & \mathbb{A}^4 &\to& \mathbb{A}^4\\
& (x,y,z,w) &\mapsto &
\left\{
\begin{array}{l}
x-(xz+yw)w\\
y+(xz+yw)z\\
z\\
w\\
\end{array}
\right. \\
\end{array}     
\]
Cette application est un contre-exemple potentiel {\`a} la conjecture des polyn{\^o}mes mod{\'e}r{\'e}s qui concerne la classification des automorphismes de l'espace affine en dimension au moins 3 \cite{Essen}.
Les applications ci-dessus, leurs conjugu{\'e}es et leur conjugu{\'e}es stable
ont un indice nef et un indice  effectif arithm{\'e}tique v{\'e}rifiant :
$\alpha(\phi, nef)\leq 2$ et $\alpha(\phi, arith.)\leq 2$ (d'apr{\`e}s les propri{\'e}t{\'e}s de stabilit{\'e} de la taille des orbites par conjugaison et conjugaison stable \cite{1art}).

\section{G{\'e}om{\'e}trie et indice ample}

\subsection{{\'E}tude de la r{\'e}solution g{\'e}n{\'e}rique}
\subsubsection{Construction}
\begin{lemm}\label{linsys}
Soit $\phi\in \Aut(\mathbb{A}^r)$ de degr{\'e} alg{\'e}brique $d >2$. Il existe $X_1$ une vari{\'e}t{\'e} projective lisse et $\psi:X_1\rightarrow \mathbb{P}^r$ un  morphisme birationnel tel que $\psi$ prolonge le morphisme $\phi:\mathbb{P}^r\setminus Z(\phi)\rightarrow \mathbb{P}^r$.

\end{lemm}

\begin{proof}
L'application $\phi:\mathbb{P}^r\rightarrow \mathbb{P}^r$  d{\'e}finie par 
$$\phi(x)=\phi(x_1,\cdots, x_{r+1})=(f_1(x),\cdots,f_{r+1}(x))$$
est une application rationnelle, lui est donc associ{\'e} un syst{\`e}me lin{\'e}aire. 
On note $V=vect(f_1,\cdots,f_{r+1})$,et  $B=\cap_{f\in V\setminus\{0\}} Z(f)$ le sous-sch{\'e}ma des points base.
Donc d'apr{\`e}s \cite{Fulton} paragraphe 4.4 nous avons le r{\'e}sultat d{\'e}sir{\'e} $\pi_1$ est l'{\'e}clatement de $\mathbb{P}^r$ le long de $B$ et $X_1$ est projectif car il peut {\^e}tre plong{\'e} dans $\mathbb{P}^r\times \mathbb{P}(V^v )$ et on identifie $\mathbb{P}(V^v )$ avec $\mathbb{P}^r$.

\noindent Nous noterons $E$ le diviseur exceptionnel associ{\'e}
\end{proof}

\begin{prop}
Soit $\phi\in \Aut(\mathbb{A}^r)$ de degr{\'e} alg{\'e}brique $d \geq 2$. Il existe $(\pi,V)$ une {\it{r{\'e}solution g{\'e}n{\'e}rique}}.
\end{prop}

\begin{proof}
Il suffit d'appliquer le lemme \ref{linsys} {\`a} $\phi$ et $\phi^{-1}$. Nous obtenons alors deux vari{\'e}tes projectives lisses $X_1$ et $X_2$, dont nous prenons le produit fibr{\'e} au dessus de $\mathbb{P}^r$, pour obtenir la vari{\'e}t{\'e} $V$.
\end{proof}

\paragraph*{Notations}Nous noterons $E$ et $F$ les diviseurs exceptionnels associ{\'e}s {\`a} cette r{\'e}solution.

\subsubsection{Propri{\'e}t{\'e}s}

\begin{prop}\label{pushforward}
Soit $\phi\in \Aut(\mathbb{A}^r)$ de degr{\'e} alg{\'e}brique $d \geq 2$. Soit $(\pi,V)$ une r{\'e}solution g{\'e}n{\'e}rique. Nous avons
$$\psi_*(H^\sharp)=\psi'_*(H^\sharp)=0.$$
De plus, si $\dim(Z(\phi))=\dim(Z(\phi^{-1}))=0$. Pour toute courbe  irr{\'e}ductible $C$ telle que $C\subset H^{\sharp}$, nous avons:
$$\psi_*(C)=\psi'_*(C)=0.$$
\end{prop}

\begin{lemm}\label{aff1}(\cite{Sibony} p.124)
Soit $\phi$ un automorphisme non-lin{\'e}aire de $\mathbb{A}^r$. Nous avons alors\,:\\
$$\phi(H\setminus Z(\phi))\subset Z(\phi^{-1})\qquad
{\mbox {et}} \qquad 
\phi^{-1}(H\setminus Z(\phi^{-1}))\subset Z(\phi).$$
\end{lemm}

\begin{proof}
Par construction, et d'apr{\`e}s le lemme \ref{aff1} nous avons :
$$\psi(H^\sharp\setminus (H^\sharp\cap E)=\phi(H\setminus Z(\phi))\subset Z(\phi^{-1}).$$
Or, $Z(\phi^{-1})$ est de codimension au moins 2, d'o{\`u} le r{\'e}sultat.
Le raisonnement est analogue dans le cas d'une courbe irr{\'e}ductible, si on suppose : $\dim(Z(\phi))=\dim(Z(\phi^{-1}))=0$.
\end{proof}

\subsection{Un r{\'e}sultat pr{\'e}liminaire}



\begin{theo}[Th{\'e}or{\`e}me de Kleiman (\cite{Lazarsfeld} th 1.4.8.)]\label{klei}
Soit $X$ une vari{\'e}t{\'e} (un sch{\'e}ma) complet. Si $D$ est un $\mathbb{R}$-diviseur nef de $X$, alors pour toute sous-vari{\'e}t{\'e} (sch{\'e}ma) irr{\'e}ductible $V\subset X$ de dimension $K$, nous avons :
$$D^K.V\geq 0.$$ 
\end{theo}

\subsection{Preuve th{\'e}or{\`e}me C}

\begin{proof}
Soit $(\pi,V)$ la r{\'e}solution g{\'e}n{\'e}rique de $\phi$.
Soit $D\subset H^{\sharp}$ une courbe irr{\'e}ductible. 
D'apr{\`e}s la proposition \ref{pushforward}, nous avons : 
$$\psi_*(C)=\psi'_*(C)=0.$$ 
Dans le cas de la dimension 2, on consid{\`e}re $H^{\sharp}$. 
Calculons {\`a} l'aide de la formule de projection:
\begin{eqnarray*}
D(\alpha,\pi,V).C&=&(\psi^*(H)+\psi'^*(H)-\alpha \pi^*(H)).C\\
&=& H.\psi_*(C)+H.\psi'_*(C) -\alpha H.\pi_*(C)=-\alpha H.\pi_*(C)\\
\end{eqnarray*}
Or, $H.\pi_*(C)>0$.
Donc en utilisant la contrapos{\'e}e du th{\'e}or{\`e}me de Kleiman \ref{klei}, nous avons :
$$\mbox{si}\quad \alpha>0\quad \mbox{alors}\quad D(\alpha,\pi,V)\quad \mbox{n'est pas un diviseur nef}.$$
Donc d'apr{\`e}s la proposition \ref{eclat} et la proposition \ref{amp=nef}, nous avons:
$$\alpha(\phi, amp)=\alpha(\phi, nef)\leq 0.$$
L'autre in\'egalit\'e est imm\'ediate \`a partir du lemme 2.1
\end{proof}

\begin{rema}\label{star}
\begin{enumerate}
\item Dans le cas de la dimension 2,
si $D(\alpha,\pi,V)$ est un $\mathbb{Q}$-diviseur
le crit{\`e}re de Nakai-Moishezon (\cite{hart} chapitre 5 th{\'e}or{\`e}me 1.10) est suffisant, on peut l'appliquer {\`a} un multiple entier de $D(\alpha,\pi,V)$.
\item Comme nous le voyons dans la preuve ci-dessus, il est possible de remplacer l'hypoth{\`e}se sur les dimensions par l'hypoth{\`e}se plus faible suivante:
\begin{center}
(H)   Il existe $C\subset H^\sharp$ une courbe irr{\'e}ductible telle que: $\psi_*(C)=\psi'_*(C)=0$.
\end{center}
Une question  naturelle se pose est-il possible d'affaiblir voir de
supprimer cette hypoth{\`e}se ?
\end{enumerate}
\end{rema}

Les automorphismes r{\'e}guliers apparaissent de mani{\`e}re naturelle en dynamique holomorphe (voir \cite{Sibony}).

\begin{defi}
Soit $\phi$ un automorphisme de $\mathbb{A}^r$ avec $r \geq 2$. L'automorphisme $\phi$
est dit {\it{r{\'e}gulier}} 
si $\deg(\phi)>1$ et \,:
$$Z(\phi)\cap Z(\phi^{-1})=\emptyset.$$
\end{defi}

\begin{lemm}(\cite{Sibony} proposition 2.3.2) \label{sib1}
Soit $\phi$ un automorphisme r{\'e}gulier de $\mathbb{A}^r$. Nous notons $l=\dim(Z(\phi^{-1}))+1$, alors :
$$\deg(\phi)^{l}=\deg(\phi^{-1})^{r-l}.$$
De plus, $\dim(Z(\phi))+\dim(Z(\phi^{-1}))=r-2$.
\end{lemm}

\begin{prop}
Soit $\phi$ un automorphisme r{\'e}gulier de $\mathbb{A}^3$. Soit $(\pi,V)$ la r{\'e}solution g{\'e}n{\'e}rique associ{\'e}e {\`a}  $\phi$.

Si $\deg(\phi)\geq\deg(\phi^{-1})$, on suppose qu'il existe une courbe irr{\'e}ductible $C\subset H^\sharp$ telle que $\psi_*(C)=0$  alors,
$$\alpha(\phi, amp)\leq 0.$$ 
\end{prop}

La m{\'e}thode est la m{\^e}me que pour d{\'e}montrer le th{\'e}or{\`e}me C, on utilise simplement en plus la remarque \ref{star} et le lemme \ref{sib1}, en effet en dimension 3 pour un automorphisme r{\'e}gulier l'un des deux lieux de non-d{\'e}finition ($Z(\phi)$  ou  $Z(\phi^{-1})$) est de dimension 0.

\begin{rema}De la m{\^e}me mani{\`e}re, nous obtenons le r{\'e}sultat suivant. Soit $\phi$ un automorphisme r{\'e}gulier de $\mathbb{A}^r$ v{\'e}rifiant $\deg(\phi)=\deg(\phi^{-1})^{r-1}$. Soit $(\pi,V)$ la r{\'e}solution g{\'e}n{\'e}rique associ{\'e}e {\`a}  $\phi$.
S'il existe une courbe irr{\'e}ductible $C\subset H^\sharp$ telle que $\psi'_*(C)=0$, alors $\alpha(\phi, amp)\leq 0.$
\end{rema}



\newpage

\vspace{4cm}
\appendix
\section{R{\'e}sultats de dynamique arithm{\'e}tique}

La dynamique arithm{\'e}tique s'int{\'e}resse aux comportements, sur $\overline{\mathbb{Q}}$, d'endomorphismes de vari{\'e}t{\'e}s, des vari{\'e}t{\'e}s lisses projectives
pour D. Northcott en 1950 \cite{Northcott}, des automorphismes de l'espace affine plus r{\'e}cemment \cite{Sil1}, \cite{Denis}, \cite{moi} et \cite{1art}.
Id{\'e}alement, il s'agit de relier la g{\'e}om{\'e}trie des applications {\`a} leur dynamique arithm{\'e}tique, tout comme dans les travaux sur le nombre de points de hauteurs born{\'e}es des vari{\'e}t{\'e}s (voir par exemple \cite{Peyre}).
Dans le th{\'e}or{\`e}me \ref{reg}, le lien entre la g{\'e}om{\'e}trie des applications et leur dynamique arithm{\'e}tique est clairement {\'e}tabli.

\subsection{R{\'e}sultats g{\'e}n{\'e}raux}
Nous donnons ici des r{\'e}sulats que nous avons utilis{\'e}s ou qui les pr{\'e}cisent. Nous renvoyons le lecteur {\`a}  \cite{moi}, et \cite{1art}.

Nous avons tout d'abord le r{\'e}sultats g{\'e}n{\'e}ral suivant:
\begin{prop}\label{hdyn}
Soit $\phi$ un automorphisme de l'espace affine
  $\mathbb{A}^r$ et $P$ un point non p{\'e}riodique de $\mathbb{A}^r(\overline{\mathbb{Q}})$.
Si $\delta(\phi)=1$, alors
$$\lim_{B\rightarrow +\infty} \frac{N(\phi,P,B)}{\log B}=+\infty.$$
\end{prop}

Nous montrons que la ''taille des orbites'' est essentiellement stable
par conjugaison par un automorphisme affine (\cite{moi} et \cite{1art}).
La notion de conjugaison stable est une g{\'e}n{\'e}ralisation de la
notion de conjugaison par un automorphisme de l'espace affine (voir
par exemple \cite{Essen}), nous montrons que la ``taille des orbites''
est {\'e}galement essentiellement stable par conjugaison stable (\cite{moi} et \cite{1art}).

Nous obtenons {\'e}galement un r{\'e}sultat plus pr{\'e}cis que la proposition \ref{hdyn} pour une certaine cat{\'e}gorie d'automorphismes de degr{\'e} dynamique 1: les automorphismes triangulaires.

\begin{theo}\label{1art}
Soit $\mu ^{\infty}$ l'ensemble des racines de l'unit{\'e}.
Soit $\phi$ une application triangulaire de $\mathbb{A}^r$,
\begin{eqnarray*}
\phi(X_{1},\dots,X_{r})&=&(a_1X_{1}+F_1(X_{2},\dots,X_{r}),a_2X_2+F_2(X_{3},\dots,X_{r}),\\
& &\quad \dots,a_{r-1}X_{r-1}+F_{r-1}(X_{r}),a_rX_r+F_r);
\end{eqnarray*}
avec   $F_i\in k[X_{i+1},\dots,X_{r}]$ pour $1\leq i\leq r-1$  et  $a_i\in
k^*$.
\newpage
Soit $P=(x_1,\dots,x_r)\in\mathbb{A}^r({\overline{\mathbb{Q}}})$ un point non
p{\'e}riodique, alors\,:
\begin{itemize}
\item s'il existe  $i$ tel que $a_i\notin \mu ^{\infty}$ et $x_i\neq
  0$ alors, $N(\phi,P,B)\gg \ll B;$
\item sinon  $\log N(\phi,P,B)\gg \ll B.$
\end{itemize}
\end{theo}

\subsection{Automorphismes r{\'e}guliers}

Une famille d'automorphismes affines se distingue naturellement
en dynamique holomorphe \cite{Sibony}, il s'agit des automorphismes r{\'e}guliers.

\begin{defi}
Soit $\phi$ un automorphisme de $\mathbb{A}^r$, $Z(\phi)$ d{\'e}signe le lieu de non-d{\'e}finition de l'application
rationnelle (d{\'e}finie sur $\mathbb{P}^r$) associ{\'e}e {\`a} $\phi$. L'automorphisme $\phi$
est dit {\it{r{\'e}gulier}} 
si $\deg(\phi)>1$ et \,:
$$Z(\phi)\cap Z(\phi^{-1})=\emptyset.$$
\end{defi}

Nous relions la dynamique
arithm{\'e}tique des automorphismes r{\'e}guliers {\`a} leur
g{\'e}om{\'e}trie.

\begin{theo}\label{reg}
Soit $\phi$ un automorphisme r{\'e}gulier de $\mathbb{A}^r$ de
degr{\'e} $d$ et $l=dim Z(\phi^{-1})+1$.
Pour tout point non p{\'e}riodique $P$ de ${\mathbb{A}}^r({\overline{\mathbb{Q}}})$, nous avons\,:
$$N(\phi,P,B)\sim \frac{r}{l}\frac{\log B}{\log d}.$$
\end{theo}
En dimension $2$, les applications de
H{\'e}non g{\'e}n{\'e}ralis{\'e}es sont des automorphismes r{\'e}guliers.

\subsection{Situation produit}

Si $\phi_i\in \Aut(\mathbb{A}^{r_i})$ pour $1\leq s$, nous d{\'e}finisons l'application
$\Phi\in \Aut(\mathbb{A}^{r_1+\cdots+r_s})$ de la fa{{\c c}}on suivante\,:
$$\Phi(x_1,\cdots,x_{r_1+\cdots+r_s})=(\phi_1(x_{1,1},\cdots,x_{r_1,1}),\cdots,
\phi_s(x_{1,s},\cdots,x_{r_s,s})).
$$
La connaissance du comportement des orbites de chaque application $\phi_i$,
nous donne le comportement des orbites de l'automorphisme $\phi$.

Nous pouvons par exemple consid{\'e}rer l'exemple suivant:
Soit $\phi_1\in \Aut(\mathbb{A}^{2})$ une application {\'e}l{\'e}mentaire. 
Soit $\phi_2$  une application de H{\'e}non g{\'e}n{\'e}ralis{\'e}e.
L'automorphisme $\Phi$ a pour
expression\,: 
$$\Phi=(\alpha x+q(y),\beta y+\gamma,p(z)-at,z).$$ 

\noindent De plus, nous avons $\delta(\Phi)=\delta(\phi_2)$.

\newpage
\noindent{\bf{Dynamique arithm{\'e}tique de $\Phi$\,:}}
\begin{itemize}
\item[$\bullet$] soit $P=(x,y,z,t)$. Si $P_1=(x,y)\notin Per(\phi_1,k)$
  et $P_2=(z,t)\in  Per(\phi_2,k)$ alors\,:
 $$N(\Phi,P,B)=N(\phi_1,P_1,B).$$ 
Suivant le type d'automorphismes
 {\'e}l{\'e}mentaires diff{\'e}rents cas se pr{\'e}sentent. 
Les tailles des orbites sont donc des
 fonctions qui asymptotiquement se comportent comme $cB$ ou $ce ^B$.
\item[$\bullet$] soit $P=(x,y,z,t)$. Si $P_1=(x,y)$
  et $P_2=(z,t)\notin  Per(\phi_2,k)$ alors\,:
$$N(\Phi,P,B)\sim 2\frac{\log B}{\log \delta(\Phi)}.$$
\end{itemize}

\subsection{Exemples}
Nous donnons ici un autre exemple d'application de degr{\'e} dynamique {\'e}gal {\`a}1: l'application d'Anick.

\begin{prop}
L'application d'Anick a pour expression:
\[
\begin{array}{lccl}
\sigma : & \mathbb{A}^4 &\to& \mathbb{A}^4\\
& (x,y,z,w) &\mapsto &
\left\{
\begin{array}{l}
x-(xz+yw)w\\
y+(xz+yw)z\\
z\\
w\\
\end{array}
\right. \\
\end{array}     
\]

Les orbites des points non p{\'e}riodiques de $\sigma$ sont
  contenues dans les quadriques d'{\'e}quations $xz+yw=c$, o{\`u} $c$ est une constante non nulle, elles ont pour taille $N(B)\sim c(x,y,z,w) e^B$.
\end{prop}

\begin{theo}\label{nagt}
L'application  de Nagata-``tordue'' est d{\'e}finie de la mani{\`e}re suivante:  

\[
\begin{array}{lccl}
\phi : & \mathbb{A}^3 &\to& \mathbb{A}^3\\
& (X,Y,Z) &\mapsto &
\left\{
\begin{array}{l}
Y-2(YZ+X^2)X-(YZ+X^2)^2Z\\
X+(YZ+X^2)Z\\
Z\\
\end{array}
\right. \\
\end{array}     
\]

Soit P un point non p{\'e}riodique, alors:
\begin{enumerate}
\item Soit P de coordonn{\'e}es (X,Y,Z) avec $Z\neq 0$,
 $$N(\phi,P,B)\sim 2\log_{4}(B).$$
 \item Soit P de coordonn{\'e}es (X,Y,0) 
$$ N(\phi,P,B)\sim 2\log_{3}(B).$$
\end{enumerate}
\end{theo}

\newpage
\backmatter

\end{document}